\documentclass{amsart}
\usepackage{amssymb}
\usepackage[english]{babel}

\newcommand{\Ps}{\mathbf{P}}
\newcommand{\Z}{\mathbf{Z}}
\newcommand{\A}{\mathbf{A}}
\newcommand{\C}{\mathbf{C}}

\newcommand{\Hu}{\mathcal{H}}
\newcommand{\str}{\mathcal{O}}
\newcommand{\ra}{\rightarrow}
\newcommand{\M}{\mathcal{M}}
\newcommand{\T}{\mathcal{T}}

\renewcommand{\phi}{\varphi}
    \newtheorem{Lem}{Lemma}[section]
    \newtheorem{Prop}[Lem]{Proposition}
    \newtheorem{Thm}[Lem]{Theorem}
    
    \newtheorem{Cor}[Lem]{Corollary}
    
    \newtheorem{Ass}[Lem]{Assumption}

   \theoremstyle{definition}
    \newtheorem{Def}[Lem]{Definition}

    \newtheorem{Rem}[Lem]{Remark}
    \DeclareMathOperator{\rank}{rank}
\DeclareMathOperator{\Sym}{Sym}
\DeclareMathOperator{\sing}{sing}
\DeclareMathOperator{\Sing}{Sing}
\DeclareMathOperator{\codim}{codim}
\DeclareMathOperator{\Aut}{Aut}
\DeclareMathOperator{\prim}{prim}
\DeclareMathOperator{\gen}{gen}
\DeclareMathOperator{\Ima}{Im}
\DeclareMathOperator{\tor}{tor}
\DeclareMathOperator{\PGL}{PGL}
\DeclareMathOperator{\NL}{NL}
\DeclareMathOperator{\MW}{MW}\DeclareMathOperator{\End}{End}
\DeclareMathOperator{\tr}{tr}

\begin{document}
\title{Higher Noether-Lefschetz loci of elliptic surfaces}
\author{Remke Kloosterman}
\address{Department of Mathematics and Computer Science, University of Gro\-ning\-en, PO Box 800, 9700 AV  Groningen, The Netherlands}
\curraddr{Institut f\"ur Al\-ge\-bra\-ische Geo\-me\-trie, Universit\"at Hannover, Wel\-fen\-gar\-ten 1, 30167 Hannover, Germany}
\email{kloosterman@math.uni-hannover.de}
\thanks{The first part of the research we report on  was done during the author's stay as EAGER pre-doc at the Turin node of EAGER. The author would like to thank Alberto Conte and Marina Marchisio for
 making this possible. The author wishes to thank Bert van Geemen for
 suggesting some literature on similar topics which influenced parts
 of this paper. The author would like to thank Joseph Steenbrink for
  conversations on the results of \cite{Stermk}. The author would like to
  thank Jaap Top and Marius van der Put for many valuable discussions on this paper. Finally the author would thank the referee for suggesting many improvements and simplifications.
This paper is based on a chapter of the author's PhD thesis \cite[Chapter 2]{proefschrift}.
}
\begin{abstract}
We calculate the dimension of the locus of Jacobian elliptic surfaces over $\Ps^1$ with  given Picard number, in the corresponding moduli space.
\end{abstract}

\date{\today}
\maketitle

\section{Introduction}
Let $\M_n$ be the coarse moduli space of Jacobian elliptic surfaces $\pi : X \ra \Ps^1$ over $\C$, 
such that the geometric genus of $X$ equals $n-1$ and $\pi$ has at least one singular fiber.
It is known that $\dim \M_n=10n-2$ (see \cite{MiMS}).  By $\rho(X)$ we denote the Picard number of $X$. It is well known that for an elliptic surface with a section we have that $2\leq \rho(X) \leq h^{1,1}=10n$.

Fix an integer $r\geq 2$, then in $\M_n$  one can study the loci
\[ \NL_r:= \{ [ \pi: X \ra \Ps^1] \in \M_n \mid \rho(X)\geq r \}. \]
We call these loci {\em higher Noether-Lefschetz loci}, in analogy with \cite{Cox}. One can show that $\NL_r$ is a countable union of Zariski closed subsets of $\M_n$ (see \cite{DelKap}). This fact can also be proven using 
the explicit description of $NS(X)$ for a Jacobian elliptic surface $\pi: X\ra \Ps^1$, due to Shioda and Tate (see Theorem~\ref{ST}).

The aim of this paper is to study the dimension of $\NL_r$.

\begin{Thm}\label{MainThmNL} Suppose $n\geq 2$. For $2\leq r \leq 10n$, we have
\[ \dim \NL_r \geq 10n-r=\dim \M_n-(r-2). \]
Moreover, we have equality when we intersect $\NL_r$ with the locus of elliptic surfaces with non-constant $j$-invariant.
\end{Thm}

The fact that the locus of elliptic surfaces with constant $j$-invariant has
dimension $6n-3$ implies
\begin{Cor}  Suppose $n\geq 2$. For $2\leq r \leq 4n+3$, we have
\[ \dim \NL_r = 10n-r. \]
\end{Cor}

Since the classes of the image of the zero-section and of a general fiber give rise to two independent classes in $NS(X)$, we have that $\NL_2=\M_n$, proving  Theorem~\ref{MainThmNL} for the case $r=2$. For $r=3$ the result was proven by Cox (\cite{Cox}). If $n=1$ then we are in the case of rational elliptic surfaces. In this case we have the well-known result $\M_1=\NL_{10}$ and $\NL_{11}=\emptyset$. If $n=2$ then we are in the case of $K3$ surfaces, and the above results follow from general results on the period map. In fact, for $K3$ surfaces we have that $\dim \NL_r=20-r$, for $2\leq r \leq 20$. We will focus on the case $n>2$.

Suppose $\pi: X \ra \Ps^1$ is an elliptic surface not birational to a product.
Denote by $\MW(\pi)$ the group of sections of $\pi$. By the Shioda-Tate formula (see Theorem~\ref{ST}) we obtain that the rank of  $\MW(\pi)$ is at most $\rho(X)-2$. {}From this and Theorem~\ref{MainThmNL} we obtain

\begin{Cor}\label{MWcor} Suppose $n\geq 2$. Let 
\[ \MW_r:= \{ [ \pi: X \ra \Ps^1] \in \M_n\mid \rank \MW(\pi)\geq
  r \}.\] Let $U:=\{[\pi: X \ra \Ps^1] \mid j(\pi) $ non-constant $ \}$. Then for $0 \leq r \leq 10n-2$ we have
\[ \dim \MW_r\cap \; U \leq 10n-r-2.\]
\end{Cor}

Cox \cite{Cox} proved that $\dim \MW_1 = 9n-1$, which is actually a stronger result than Corollary~\ref{MWcor} for the special case $r=1$.

The proof of Theorem~\ref{MainThmNL} consists of two parts. In the first part we construct elliptic surfaces with high Picard number. This is done by constructing large families of elliptic surfaces such that the singular fibers have many components. To calculate the dimension of the locus of this type of families, we study the ramification of the $j$-map, and calculate the dimension of several Hurwitz spaces. This yields $\dim \NL_r\geq 10n-r$.

The second part consists of proving that $\dim \NL_r \cap \; U \leq 10n-r$. 
We choose a strategy similar to what  M.L. Green \cite{GreenF} uses in order to identify the components of maximal dimension in the Noether-Lefschetz locus in the case of surfaces of degree $d$ in $\Ps^3$. In order to apply this strategy we consider an elliptic surface over $\Ps^1$ with a section  as a surface $Y$ in the weighted projective space $\Ps(1,1,2n,3n)$ with $n=p_g(X)+1$. To obtain $Y$, we need to contract the zero-section and all fiber components not intersecting the zero-section.
Then we use Griffiths' and Steenbrink's identification of the Hodge filtration on $H^2(Y,\C)$ with graded  pieces of the Jacobi-ring of $Y$ (the coordinate ring of $\Ps$ modulo the ideal generated by the partials of the defining polynomial of $Y$). Using some results from commutative algebra we can calculate an upper bound for  the dimension of $\NL_r\cap \; U$.

We would like to point out an interesting detail: the classical Griffiths-Steenbrink identification holds under the assumption that $Y$ is smooth outside the singular locus of the weighted projective  space. In  our case it might be that $Y$ has finitely many rational double points  outside the singular locus of the weighted projective space. Recently, Steenbrink (\cite{Stermk}) obtained a satisfactory identification in the case that $Y$ has ``mild'' singularities.

When we consider elliptic surfaces with constant $j$-invariant 0 or 1728, the theory becomes slightly more complicated: 
several families of elliptic surfaces over $\Ps^1$ with a section and constant $j$-invariant 0 or 1728 and generically Picard number $\rho$ have  codimension one subfamily of surfaces with Picard number $\rho+2$.
It turns out that for certain values of $r\geq 8n$, such families
prevent us from proving the equality $\dim \NL_r=10n-r$. These surfaces have other strange properties. For the same reason as above, we can produce examples not satisfying several Torelli type theorems (see \cite[Theorem 4.8]{Ext}). These surfaces are also the elliptic surfaces with larger Kuranishi families than generic elliptic surfaces. Actually, the difference between the dimension of $\NL_r$ and $10n-r$ equals the difference between the dimension of the Kuranishi family of a generic elliptic surface in $\NL_r$ and the 
dimension of the Kuranishi family of a generic elliptic surface $X$ over $\Ps^1$, with $p_g(X)=n-1$ and admitting a section.

The organization of this paper is as follows. In Section~\ref{PrelNL} we recall several standard facts on elliptic surfaces. 
In Section~\ref{Hur} we calculate the dimension of certain Hurwitz spaces. In Section~\ref{config} we study configurations of singular fibers. In Section~\ref{under} we use the results of the previous two sections to identify several components in $\NL_r$ of dimension $10n-r$.
In Section~\ref{constj} we study the locus in $\M_n$ of elliptic surfaces with constant $j$-invariant. In Section~\ref{jnul} we study elliptic surfaces with the `special' $j$-invariants 0 and 1728. We use these surfaces to identify components $L$ of $\NL_r$ such that $\dim L \gg 10n-r$.
In Section~\ref{upper} we prove that the identified components are of maximum  dimension in $\NL_r$. This is done by applying a modified version of the Griffiths-Steenbrink identification of the Hodge structure of hypersurfaces with several graded pieces of the Jacobi-ring.
In Section~\ref{concl} some remarks are made and some open questions
are raised.
\section{Definitions \& Notation}
\label{PrelNL}

\begin{Ass} In the sequel we work over the field of complex numbers.

By a curve we mean a non-singular projective complex connected curve.

By a surface we mean a projective complex surface. Moreover, if not specified otherwise, we assume it to be smooth.\end{Ass}

\begin{Def} A \emph{(Jacobian) elliptic surface (over $\Ps^1$)} is a morphism $\pi: X\ra \Ps^1$ together with a section $\sigma_0: \Ps^1\ra X$ to $\pi$, with $X$ a relatively minimal surface and such that almost all fibers are  irreducible genus 1 curves.

  We denote by $j(\pi): \Ps^1 \rightarrow \Ps^1$ the rational function such that $j(\pi)(P)$ equals the $j$-invariant of $\pi^{-1}(P)$, whenever $\pi^{-1}(P)$ is non-singular.

  The set of sections of $\pi$ is an abelian group, with $\sigma_0$ as  identity element. Denote this group by $\MW(\pi)$.

\end{Def}

\begin{Def} Let $\pi: X \rightarrow \Ps^1$ be an elliptic surface. Let $P$ be a point of $\Ps^1$. Define $v_P(\Delta_P)$ as the valuation at $P$ of the minimal discriminant of the Weierstrass model, which equals the topological Euler characteristic of $\pi^{-1}(P)$. 

The \emph{degree of the elliptic surface $\pi: X\ra \Ps^1$} is the degree of the line bundle $[R^1\pi_*\str_X]^{-1}$.
\end{Def}

\begin{Prop}\label{Noether} Let $\pi: X \rightarrow \Ps^1$ be an elliptic surface of degree $n$. Then
\[ \sum_{P\in \Ps^1} v_P(\Delta_P) = 12n. \]
\end{Prop}

\begin{proof} This follows from Noether's formula (see \cite[p. 20]{BPV}). The precise reasoning can be found in \cite[Section III.4]{MiES}.
\end{proof}

\begin{Def} Let $X$ be a surface, let $C$ and $C_1$ be curves. Let $\varphi: X \rightarrow C$ and $f: C_1 \rightarrow C$ be two morphisms. Then we denote by $\widetilde{X \times_C C_1}$ the smooth, relatively minimal model of the ordinary fiber product of $X$ and $C_1$.

We call the induced morphism $\widetilde{X\times_C C_1}\ra C_1$ the base-change of $\varphi$ by $f$.
\end{Def}

\begin{Def} For a surface $X$. Let $NS(X)$ denote the group of divisor on $X$ modulo algebraic equivalence. We call $NS(X)$ \emph{the N\'eron-Severi group} of $X$.

Let  $\rho(X)$ denote the rank of the
  N\'eron-Severi group of $X$. We call $\rho(X)$ the {\em Picard number}.\end{Def}

One can show that if the degree of $\pi: X \ra \Ps^1$ is positive then it equals $p_g(X)+1$. (See \cite[Lemma III.4.2]{MiES}.) If the degree is not positive then it is zero and $\pi$ is the projection $E\times \Ps^1\ra\Ps^1$, for some elliptic curve $E$.

In the sequel we need a description of the N\'eron-Severi group. 
There is a simple description of the N\'eron-Severi group for elliptic surfaces. Recall the following theorem.
\begin{Thm}[{Shioda-Tate \cite[Theorem 1.3 \& Corollary 5.3]{Sd}}]\label{ST}  Let $\pi:X\rightarrow \Ps^1$ be an elliptic surface, with at least one singular fiber. Then the N\'eron-Severi group of $X$ is generated by the classes of $\sigma_0(\Ps^1)$, a fiber, the components of the singular fibers not intersecting $\sigma_0(\Ps^1)$, and the generators of the Mordell-Weil group. Moreover, let $S$ be the set of points $P$ such that $\pi^{-1}(P)$ is singular. If we denote by $m(P)$  the number of irreducible components of $\pi^{-1}(P)$, then
  \[ \rho(X) := \rank(NS(X))=2+ \sum_{P \in S} (m(P)-1)+\rank(\MW(\pi)) \]
\end{Thm}

\begin{Def} Suppose $\pi: X \ra \Ps^1$ is an elliptic surface. Denote by $T(\pi)$ the  subgroup of the N\'eron-Severi group of $X$ generated by the classes of the fiber, $\sigma_0(\Ps^1)$ and the components of the singular fibers not intersecting $\sigma_0(\Ps^1)$. Let $\rho_{\tr}(\pi):=\rank T(\pi)$. We call $T(\pi)$ the {\em trivial part} of the N\'eron-Severi group of $X$.
\end{Def}

\begin{Rem}In Lemma~\ref{redlem} we give an alternative description for the trivial part of the N\'eron-Severi group.\end{Rem}

\begin{Def} Suppose $n\geq 1$ is an integer. We denote by $\M_n$ the moduli
  space of Jacobian elliptic surfaces $\pi: X \ra \Ps^1$ of degree $n$.\end{Def}

One can show that $\M_n$ is a quasi-projective variety of dimension $10n-2$. The moduli space $\M_n$ is constructed using Geometric Invariant Theory, and one can show that for $n>2$ all elliptic surfaces are stable. This implies that all isomorphism classes of elliptic surfaces over $\Ps^1$ of degree $n>1$ yield points in $\M_n$. For more information on these moduli spaces we refer to \cite{MiMS}.
\section{Dimension of Hurwitz Spaces}\label{Hur}
In this section we calculate the dimension of several Hurwitz spaces. We expect that all the results in this section are already known to the experts. Unfortunately, we could not find an exposition on this subject in the literature which would be sufficient for the application in the sequel. Many of the ideas used in this section are also present in \cite{MiRa}.

\begin{Def} Let $C_1$ and $C_2$ be curves. Two morphisms $\varphi_i:C_i \ra \Ps^1$ are called isomorphic, if there exists an isomorphism $\psi: C_1 \ra C_2$ such that $\varphi_1=\varphi_2\circ \psi$.\end{Def}

\begin{Def} Let $m>2$ be an integer. Fix $m$ distinct points $P_i\in \Ps^1$. Let $\Hu(\{e_{i,j}\}_{i,j})$ be the Hurwitz space (coarse moduli space) of isomorphism classes of semi-stable morphisms $\varphi: \Ps^1 \rightarrow \Ps^1$ of degree $d$ such that $\varphi^*P_i=\sum_k e_{i,j} Q_j$, with $Q_{j'} \neq Q_j$ for $j'\neq j$. (If $\varphi$ satisfies this condition, we will say that $\varphi$ has \emph{ramification indices $e_{i,j}$  over the $P_i$.})\end{Def}

\begin{Rem} Note that morphisms corresponding to points  of the Hurwitz space $\Hu(\{e_{i,j}\}_{i,j})$ might be ramified outside the $P_i$.\end{Rem}

In the following remark we indicate where the notion `semi-stable' comes from.
\begin{Rem}\label{ss} 
All morphisms of degree $d$ can be parameterized by an open set in $\Ps^{2d+1}$ by sending a point $[x_0:x_1:\dots:x_{2d+1}]$ to the morphism induced by the function $t\mapsto (x_0+x_1t+\dots +x_d t^d )/(x_{d+1}+x_{d+2} t+\dots + x_{2d+1}t^d)$. It is not hard to write down a finite set of equations and inequalities in the $x_i$, such that every solution corresponds to a morphism with the required ramification behavior over the $P_i$. To obtain $\Hu(\{e_{i,j}\}_{i,j})$, one needs to divide out by the action of the reductive group $\Aut(\Ps^1)=\PGL_2$. Geometric Invariant Theory ensures the existence of a `good' quotient, possibly after restricting to the smaller open subset of of so-called \emph{semi-stable} elements for the action of $\PGL_2$. (For more information of construction of moduli spaces of this type and a precise definition of semi-stability, see for example \cite{dolinv}.)

We do not describe which morphisms are semi-stable. Indeed, since we are only interested in the dimension of $\Hu(\{e_{i,j}\}_{i,j})$, it is enough for our purposes to work with a dense open subset of $\Hu(\{e_{i,j}\}_{i,j})$. In the case that a morphism with the prescribed ramification exists, we construct  for some $k\geq 0$ a finite \'etale covering of a Zariski open in $\Sym^k(\Ps^1)$ which parameterizes isomorphism classes of morphisms with simple ramification outside the $P_i$. It turns out that $\Hu(\{e_{i,j}\}_{i,j})$ is a partial compactification of this space.
\end{Rem}

\begin{Rem} If $m>3$, then $\Hu(\{e_{i,j}\}_{i,j})$ depends on the points $P_i$. We will prove that its dimension is independent of the choice of the $P_i$. \end{Rem}

Recall the following special form of the Riemann existence theorem.
\begin{Prop}\label{Riemann} Fix $m\geq 2$ points $R_i\in \Ps^1$. Fix a positive integer $d$. Fix partitions of $d$ of the form $d=\sum_j^{k_i} e_{i,j}$, for $i=1,\dots, m$. Let $q=\sum k_i$. Assume that $q=(m-2)d+2$. Then we have a correspondence (the so-called monodromy representation) between
\begin{itemize}
\item Isomorphism classes of morphisms $\varphi:\Ps^1\ra\Ps^1$ with ramification indices $e_{i,j}$ over the $R_i$ and unramified elsewhere.

\item Congruence classes of transitive subgroups of $S_d$ (the symmetric group on $d$ letters) generated by $\sigma_i, i=1,\dots,m$, such that the lengths of the cycles of $\sigma_i$ are the $e_{i,j},j=1,\dots k_i$ and $\prod \sigma_i=1$.
\end{itemize}
\end{Prop}

\begin{proof} In \cite[Corollary 4.10]{MiRS} the above equivalence is proven, except that they consider all morphisms $C \ra \Ps^1$ with given ramification indices, and they do not assume $q=(m-2)d+2$. Hence we need to show that $g(C)=0$ is equivalent to $q=(m-2)d+2$: The condition $g(C)=0$ is equivalent by the Hurwitz' formula (see \cite[Corollary IV.2.4]{Har})  to
\[ -2=-2d+\sum_{Q\in \Ps^1} e_Q(\varphi)-1,\]
where $e_Q(\varphi)$ is the ramification index of $\varphi$ at $Q$. Since $\sum e_Q(\varphi)-1=md-q$, it follows that the displayed formula  is equivalent to $q=(m-2)d+2$, which yields the proof.
\end{proof}

\begin{Def} Let $\varphi: C \ra C'$ be a non-constant morphism of curves. Let $S$ be a set of points on $C'$. We say that $\varphi$ has {\em simple ramification outside $S$} if for every point $Q\not \in S$ we have $\# \varphi^{-1}(Q) \geq \deg(\varphi)-1$.\end{Def}

The following Lemma is called `deformation of a function' by Miranda (see \cite[Section 3]{MiRS}).
\begin{Lem} \label{ramred}Let $\varphi: C \ra \Ps^1$ be a non-constant morphism. Let $P\in \Ps^1$ be a critical point and $Q_i, i=1,\dots, s$ be the points in $\varphi^{-1}(P)$. Suppose $e_{Q_1}>1$. Then for every integer $k$ with $1\leq k \leq e_{Q_1}$ there exists a non-constant morphism $\varphi': C' \ra \Ps^1$ such that the ramification behavior  of $\varphi'$ and $\varphi$ coincide at every point $T\in \Ps^1$ except at the point $P$ and one other point $P'\in \Ps^1$. At these points we have that
\begin{itemize}
\item the morphism $\varphi$ is unramified at $P'$;
\item over $P'$ the morphism $\varphi'$ has simple ramification;
\item  we can write  $\varphi'^{-1}(P)=\{Q'_0,\dots Q'_s\}$ such that $e_{Q'_0}=k$, $e_{Q'_1}=e_{Q'_1}-k$, and $e_{Q'_j}=e_{Q_j}$ for $j=2,\dots s$.
\end{itemize}
In particular, $g(C)=g(C')$.
\end{Lem}

\begin{proof}
Proposition~\ref{Riemann} implies that we can associate with each critical point $R$ of $\varphi$ an element $\sigma_R$ of $S_d$ such that the lengths of the cycles of $\sigma_R$ coincide with the ramification indices over $R$, the subgroup generated by the $\sigma_R$ acts transitively on $\{1,\dots,d\}$ and $\prod \sigma_R =1$. 

Without loss of generality we  may assume that $\sigma_{P}$ contains the cycle $(1 \;2\; \dots\; e_{Q_1})$. Let $\tau_{P}$ be the cycle obtained from $\sigma_{P}$ by replacing $(1\; 2 \;\dots \;e_{Q_1})$ with $(1 \; 2 \;\dots \;k)(k+1\;\dots\; e_{Q_1})$. Let $\tau_{P'}$ be the cycle $(k\; e_{Q_1})$. For all critical values $R$ of $\varphi$ different from $P$ set $\tau_R=\sigma_R$. Then one easily shows that the subgroup generated by the $\tau_R$ is transitive and $\prod \tau_R=1$.
Proposition~\ref{Riemann} implies that we can associate a morphism $\varphi': C' \ra \Ps^1$ with the subgroup generated by the $\tau_R$. The same Proposition implies  that $\varphi'$ has the desired ramification behavior.

The statement on the genus of $C$ and $C'$ follows directly form the Hurwitz' formula \cite[Corollary IV.2.4]{Har}.\end{proof}

\begin{Lem}\label{simpleram}
Fix  integers $d\geq 2$ and $m\geq 2$ and fix $m$ partitions $d$ of the from $d=\sum_{j=1}^{k_i} e_{i,j}$, $i=1,2,\ldots, m$. Let $q=\sum k_i$. Fix $m$ points $P_1,\dots P_m \in \Ps^1$.

Assume there is a morphism $\varphi':\Ps^1\ra\Ps^1$ with ramification indices $e_{i,j}$ over the $P_i$. 
Then there is a morphism $\varphi$ such that $\varphi$ has simple ramification outside the $P_i$, and ramification indices $e_{i,j}$ over the $P_i$. 
\end{Lem}

\begin{proof} 
Suppose $\varphi'$  has non-simple ramification outside the $P_i$. 
Applying  Lemma~\ref{ramred} several times yields that there is a morphism  with the same ramification indices over the $P_i$ as $\varphi'$ and simple ramification elsewhere. 
\end{proof}

The above results enable us to calculate the dimension of the Hurwitz space.
\begin{Prop} \label{Hurwitzdim}
Fix  integers $d\geq 2$ and $m\geq 2$ and fix $m$ partitions $d$ of the from $d=\sum_{j=1}^{k_i} e_{i,j}$, $i=1,2,\ldots, m$. Let $q=\sum k_i$. Fix $m$ points $P_1,\dots P_d$.

Let $q=\sum k_i$.
The dimension of $\Hu(\{e_{i,j}\}_{i,j})$ is 
$q-(m-2)d-2$, provided that there exists a morphism $\varphi:\Ps^1\ra \Ps^1 $ with ramification indices $e_{i,j}$ over the $P_i$. 
\end{Prop}

\begin{proof} We prove this theorem by induction on $q-(m-2)d-2$. If $q-(m-2)d-2<0$ then the Hurwitz' formula  \cite[Corollary IV.2.4]{Har} implies that $\Hu(\{e_{i,j}\})_{i,j}$ is empty.
If $q-(m-2)d-2=0$ then the Hurwitz' formula implies that  $\varphi$ is not ramified outside the $P_i$. Proposition~\ref{Riemann} implies that $\Hu(\{e_{i,j}\}_{i,j})$ is isomorphic to a  finite collection of congruence classes of subgroups of $S_d$. This yields this case.

We prove now the general case. Let $\Hu'(\{e_{i,j}\}_{i,j})$ be the moduli space of all morphisms  $\psi: \Ps^1\ra \Ps^1$ such that the ramification indices over the $P_i$ are $e_{i,j}$ and simple ramification elsewhere. Let $\Delta$ be the 
  complement of $\Hu'(\{e_{i,j}\}_{i,j})$ in $\Hu(\{e_{i,j}\}_{i,j})$. We will use the induction hypothesis to prove that $\dim \Delta \leq q-(m-2)d-3$. We start by proving that $\dim \Hu'(\{e_{i,j}\}_{i,j})_{i,j}=q-(m-2)d-2$.
 
 Let $S$ be the collection of congruence classes of subgroups of $S_d$, associated with morphisms with the same ramification behavior as $\psi$ (cf. Proposition~\ref{Riemann}). Note that $\psi$ is ramified at $q-(m-2)d-2$ points outside the $P_i$. 

Let $U\subset \Sym^{q-(m-2)d-2} \Ps^1$ be the set of points $Q_1+\dots+ Q_{q-(m-2)d-2}$ such that for all $i,j$ we have $Q_i\neq Q_j$ and $Q_i\neq P_j$. Then Proposition~\ref{Riemann} implies that $\Hu'(\{e_{i,j}\}_{i,j})=S\times U$. In particular $\dim \Hu'(\{e_{i,j}\}_{i,j})=q-(m-2)d-2$. 

It remains to  bound the dimension of $\Delta$. Note that $\Delta$ corresponds to morphisms with some non-simple ramification outside the $P_i$. 
Fix a morphism $\varphi$ corresponding to a point in $\Delta$. Then there are $t>0$ points $Q_k\in \Ps^1$ different from the $P_i$ over which $\varphi$ has non-simple ramification.
Let $\Hu''$ be the Hurwitz space of morphisms such that the ramification indices over the $P_i$ and $Q_k$ are the same as $\varphi$.  
Then by induction we have $\dim \Hu''\leq q-(m-2)d-2-2t $. Letting the $Q_k$ move on $\Ps^1$  yields a variety of dimension $q-(m-2)d-2-t$. Since there are finitely many possibilities for the ramification indices over the points different from the $P_i$  it follows that $\Delta$  has dimension at most $q-(m-2)d-3$.

This finishes the proof.
\end{proof}

\begin{Cor} The dimension of $\Hu(\{e_{i,j}\}_{i,j})$ is independent of the choice of the $P_i$.\end{Cor} 

Suppose we know the ramification indices modulo some integer $N_i$.
The following corollary tells us that if  for one choice of the ramification indices, the associated Hurwitz space is non-empty, then the same holds for the Hurwitz space associated with the minimal choice of ramification indices. In particular,  the  Hurwitz space associated with that particular choice is the largest one. 
\begin{Cor}\label{ModN} Let $m,d$ be positive integers. Fix $m$ integers $N_i$ such that $N_i\leq d$. Let $a_{i,j}$ be integers such that $1\leq a_{i,j} < N_i$, and $r_iN_i+\sum_{j=1}^{s_i} a_{i,j} =d$, with $r_i$ a non-negative integer. Fix $m$ points $P_i$  on $\Ps^1$.

For all $i=1,\dots m$, set
\[ e_{i,j} = \left\{ \begin{array}{ll} a_{i,j} & 1 \leq j \leq s_i,\\ N_i & s_i+1\leq j \leq s_i+r_i. \end{array} \right.\]

Suppose there exist $m$ partitions $d=\sum_j^{s'_i} e'_{i,j}$ such that $s'_i\leq s_i$ and $e'_{i,j} \equiv e_{i,j} \bmod N_i$ if $1 \leq j\leq s'_i$. 
Then $\dim \Hu(\{e'_{i,j}\}_{i,j}) \leq \dim \Hu(\{e_{i,j}\}_{i,j})$ holds.
\end{Cor}

\begin{proof}Applying Lemma~\ref{ramred} sufficiently many times yields that
if the Hurwitz' scheme $\Hu(\{e'_{i,j}\}_{i,j})$ is non-emp\-ty then  $\Hu(\{e_{i,j}\}_{i,j})$ is non-emp\-ty. Now apply Proposition~\ref{Hurwitzdim}.
\end{proof}

\section{Configuration of singular fibers} \label{config}
Fix some $n\geq 2$. In this section we calculate the dimension of the locus in $\M_n$ corresponding to elliptic surfaces with a fixed configuration of singular fibers, containing a fiber of type $I_\nu$ or $I^*_\nu$, with $\nu>0$. For more on this see also \cite[Lectures V and X]{MiES} and~\cite[Sections 5 and 6]{Ext}. 

We start with introducing the notion of `twisting' which is well-known in the theory of elliptic curves.

Given an elliptic surface  $\pi: X \ra \Ps^1$, we can associate an
 elliptic curve in $\Ps^2_{\C(\Ps^1)}$  corresponding to the generic fiber
 of $\pi$. This induces a bijection between isomorphism classes of Jacobian
 elliptic surfaces and elliptic curves over $\C(\Ps^1)$.

Two elliptic curves $E_1$ and $E_2$ are isomorphic over $\C(\Ps^1)$ if and
only if $j(E_1)=j(E_2)$ and the quotient of the minimal discriminants of $E_1/\C(\Ps^1)$ and $E_2/\C(\Ps^1)$ is a 12-th power (in $\C(\Ps^1)^*$).

Assume that $E_1$, $E_2$ are elliptic curves over $\C(\Ps^1)$ with
$j(E_1)=j(E_2)\neq 0,1728$. (For example,  elliptic surfaces with a fiber of type $I_\nu$ or $I^*_{\nu}$ have non-constant $j$-invariant, hence they satisfy  this assumption.) One easily shows that
$\Delta(E_1)/\Delta(E_2)$ equals $u^6$, with $u\in \C(\Ps^1)^*$. Hence $E_1$ and $E_2$ are isomorphic over $\C(\Ps^1)(\sqrt{u})$. We call $E_2$ the twist
of $E_1$ by $u$, denoted by $E_1^{(u)}$. Actually, we are not interested in the function $u$, but in the places at which the valuation of $u$ is odd.

\begin{Def} Let $\pi:X \ra \Ps^1$ be a Jacobian elliptic surface. Fix $2N$ points $P_i \in \Ps^1$. Let $E/\C(\Ps^1)$ be 
the generic fiber of $\pi$.

A Jacobian elliptic surface $\pi': X' \ra \Ps^1$ is called a \emph{(quadratic) twist} of $\pi$ by $(P_1,\ldots,P_{2N})$ if the generic fiber of $\pi'$ is isomorphic to $E^{(f)}$, where $E^{(f)}$ denotes the quadratic twist of $E$ by $f$ in the above mentioned sense
and $f\in\C(\Ps^1)$ is a function such that  $v_{P_i}(f) \equiv 1 \bmod 2$ and $v_Q(f) \equiv 0 \bmod 2$ for all $Q\not \in \{P_i\}$.
\end{Def}
The existence of a twist of $\pi$ by $(P_1,\ldots,P_{2N})$ is immediate. One can show that the function $f$ mentioned in the above definition is unique up to squares, implying that a twist  is unique up to an isomorphism of the fibration $\pi'$. This property depends on the choice of our base curve. If we replaced $\Ps^1$ by an arbitrarily base curve $C$ we would have $2^{2g(C)}$ twists by a fixed set of points.

If $P$ is one of the $2N$ distinguished points, then the fiber of $P$ changes in the following way (see \cite[V.4]{MiES}).
\[
I_\nu \leftrightarrow I^*_\nu \;(\nu \geq 0) \;\;\; \;\;
II \leftrightarrow    IV^* \;\;\;\;\;
III  \leftrightarrow  III^* \;\;\;\;\;
IV  \leftrightarrow  II^*
\]

\begin{Def} A {\em configuration of singular fibers} 
is a formal sum $C$ of Kodaira types of singular fibers, with non-negative integer coefficients.

Let $i_\nu(C)$ denote the coefficient of $I_\nu$ in $C$. Define $ii(C)$, $iii(C)$, $iv(C)$, $iv^*(C)$, $iii^*(C)$, $ii^*(C)$ and $i_\nu^*(C)$ similarly.

A configuration $C$ satisfies {\em Noether's condition} if
\[ \sum_{\nu>0} \nu i_\nu+\sum_{\nu\geq 0} (\nu+6) i_\nu^*+2ii+3iii+4iv+8iv^*+9iii^*+10ii^*=12n(C) \]
with $n(C)$ a positive integer.

A \emph{multiplicative fiber} is a fiber of type $I_\nu, \nu>0$, an \emph{additive fiber} is a singular fiber not of type $I_\nu$. 
\end{Def}

The Kodaira types of singular fibers can be found at many places in the literature, e.g., \cite{BPV} or \cite{MiES}.

With an elliptic surface $\pi:X \rightarrow \Ps^1$ corresponding to a point in $\M_n$ we can associate its (total) configuration of singular fibers $C(\pi)$.
Then $C(\pi)$  satisfies Noether's condition, with $n(C(\pi))=n$ (this follows from Lemma~\ref{Noether}).

\begin{Ass}For the rest of the section, let $C$ be a configuration of singular fibers satisfying Noether's condition with $n(C)=n$ and containing at least one fiber of type $I_\nu$ or $I^*_\nu$, with $\nu>0$.\end{Ass}

\begin{Lem}\label{rambeh}
Suppose there exists a morphism $\varphi:\Ps^1 \rightarrow \Ps^1$, such that the ramification indices are as follows:
\begin{itemize}
\item  Above 0
\begin{itemize}
\item there are precisely $ii(C)+iv^*(C)$ points with ramification indices 1 modulo
3 and
\item there are precisely $iv(C)+ii^*(C)$ points with ramification indices 2 modulo 3.\end{itemize}

\item Above 1728 there are precisely $iii(C)+iii^*(C)$ points with ramification indices 1 modulo 2.

\item Above $\infty$ there are for every $\nu>0$ precisely $i_\nu(C)+i_\nu^*(C)$ points with ramification index $\nu$.
\end{itemize}
Then there exists an elliptic surface such that $C(\pi)=C$.

Conversely, if there exists an elliptic surface with $C(\pi)=C$, then $j(\pi)$
satisfies the above mentioned conditions.
\end{Lem}

\begin{proof} The last part of the statement follows from \cite[Lemma IV.4.1]{MiES}.

To prove the existence of $\pi$: Let $\pi_1: X_1\rightarrow \Ps^1$ be an elliptic surface with $j(\pi_1)=t$, with $t$ a local coordinate on the base curve $\Ps^1$. (For example one can take the elliptic surface associated to $y^2+xy=x^3-36/(t-1728)x-1/(t-1728)$.)

 Let $\pi_2:X_2\ra \Ps^1$ be the base-change  of $\pi_1$ 
by $\varphi$. 
Then it follows from \cite[Lemma IV.4.1]{MiES} that $i_\nu(C(\pi_2))+i_\nu^*(C(\pi_2))=i_\nu(C)+i_\nu^*(C)$, for $\nu>0$, and $ii(C(\pi_2))+iv^*(C(\pi_2))=ii(C)+iv^*(C)$, and that similar relations hold for $(iii,iii^*)$ and $(iv,ii^*)$.

It is easy to see  that there exists a twist $\pi_3$ of $\pi_2$ such that $C(\pi_3)-C=\epsilon I_0^*$, with $\epsilon\in \{0,1\}$. Since both configurations satisfy Noether's condition, it follows that $\epsilon=0$. Hence $\pi_3$ is the desired Jacobian elliptic surface.
\end{proof}

\begin{Lem}\label{io}
Assume that there exists an elliptic surface $\pi': X' \ra \Ps^1$ with $C(\pi')=C$. Then
\[ \dim \{ [\pi: X\ra \Ps^1] \in \M_{n} \mid j(\pi)=j(\pi'), C(\pi)=C \} = i_0^*(C). \]
\end{Lem}

\begin{proof}
Fix one $\pi_0: X_0 \ra \Ps^1$, with $C(\pi_0)=C$ and $j(\pi_0)=j(\pi')$.

Fix $i_0^*(C)$ points $P_i \in \Ps^1$, none of them in $j(\pi_0)^{-1}(\{0,1728,\infty\})$, such that $\pi^{-1}(P_i)$ is smooth for all $i$. Let $Q_i$ be the points over which the fiber of $\pi$ is of type $I_0^*$. Then twisting $\pi$ by the points $\{ P_i, Q_i \}$ gives an elliptic surface $\pi$ with $j(\pi)=j(\pi')$ and $C(\pi)=C$ (see Lemma~\ref{rambeh}). If two such twists are isomorphic then the set of points $\{ P_i \}$ are the same.  So
\[ \dim \{ [\pi: X\ra \Ps^1] \in \M_{n} \mid j( \pi)=j(\pi'), C(\pi)=C \} \geq  i_0^*(C). \]

As we remarked above a twist by a fixed set of points is unique. {}From this it follows that the number of twists $\pi''$ of $\pi'$ such that $\Sing(\pi')=\Sing(\pi'')$ and $C(\pi')=C(\pi'')$ is finite, where $\Sing(\psi)=\{ P\in \Ps^1 \mid \psi^{-1}(P) \mbox{ is singular}\}$.
Since all singular fibers not of type $I_0^*$ lie in \[j^{-1}(0,1728,\infty)\] it follows that 
\[ \dim \{ [\pi: X\ra \Ps^1] \in \M_{n(C)} \mid j( \pi)=j(\pi'), C(\pi)=C \} \leq  i_0^*(C). \]
Combining both bounds yields the lemma.
\end{proof}

\begin{Lem}~\label{dimLoc} Assume that there exists an elliptic surface $\pi': X' \ra \Ps^1$ with $C(\pi')=C$. Then the  locus $L(C)$ in $\M_{n(C)}$ corresponding to all elliptic surfaces with $C(\pi)=C(\pi')$ is constructible and has dimension
\[ \# \{\mbox{singular fibers}\} + \# \{ \mbox{fibers of type } II^*, III^*, IV^*, I_\nu^*\} -2n(C)-2 . \]\end{Lem}

\begin{proof} {}From the above lemmas it follows that $L(C)$ is a finite union of Zariski open subsets $U_i$ in $(\Ps^1)^{i_0^*(C)}$-bundles over some $\Hu(\{e_{i,j}\}_{i,j})$. This proves the constructibility of $L(C)$.

Let $\pi: X \ra \Ps^1$ be an elliptic surface corresponding to a point in $L(C)$.  {}From Lemma~\ref{rambeh} we obtain that the degree $d$ of $j(\pi)$ equals $\sum_\nu \nu (i_\nu+i_\nu^*)$. 

Similarly, we obtain congruence relations for the ramification indices of $j(\pi)$ over 0 and 1728. We would like to calculate the maximum of the dimensions of all Hurwitz spaces associated with different solutions of these congruence relations. {}From Corollary~\ref{ModN} it follows that we only have to consider the solution with the lowest ramification indices, i.e., the solution such that the number of points over 0, 1728 and $\infty$ is maximal.
One easily shows that over $\infty$ there are $\sum_{\nu>0} i_{\nu}+i^*_{\nu}$ points, over 0 there are $ii+iv^*$ points with ramification index 1, $ii^*+iv$ points with ramification index 2 and  $(d-ii-iv^*-2iv-2ii^*)/3$ points with ramification index 3. Over 1728 we obtain that there are $iii+iii^*$ points with ramification index 1 and $(d-iii-iii^*)/2$ points with ramification index 2. This implies that the $q$ mentioned in Proposition~\ref{Hurwitzdim} equals
\begin{Small}\[ ii+iv^* + iv+ii^* + (d-ii-iv^*-2iv-2ii^*)/3 + iii +iii^* + (d-iii-iii^*)/2 +\sum_{\nu>0}( i_\nu+i_\nu^*). \]\end{Small}

Since the number $m$ of points with prescribed ramification for $j$-invariant is  $3$, it follows from Corollary~\ref{ModN} that the union of Hurwitz spaces corresponding to $j$-invariants giving rise to elliptic surfaces in $L(C)$ has dimension $q-d-2$, hence Lemma~\ref{io} implies that
\[ \dim L(C)=q-d-2+i_0^*(C). \]

A simple calculating using Noether's condition yields that 
\begin{Small}
\begin{eqnarray*} q-d-2 &=&\frac{  8ii+8iv^* + 6iii+6iii^*+4iv + 4ii^* + 12 \sum_{\nu>0} (i_\nu+i_\nu^*) - 2d-24}{12} \\
&=& \sum_{\nu>0} (i_\nu+2i_\nu^*)+ii+iii+iv+i_0^*+2iv^*+2iii^*+2ii^*  -2 -2n.
\end{eqnarray*}
\end{Small}

This implies the lemma.
\end{proof}

\begin{Prop}\label{dimlem} Let $C$ be a configuration of singular fibers,  containing at least one $I_\nu$ or $I_\nu^*$-fiber ($\nu>0$) and such that there exists an elliptic surface $\pi': X' \ra \Ps^1$ with $C(\pi')=C$. Then
the dimension of $\{[\pi: X\rightarrow \Ps^1] \in \M_{n} \mid C(\pi)=C\}$ equals  \[ 10n-\rho_{\tr}(\pi)-\# \{\mbox{fibers of type } II, III \mbox{ or } IV \}.\]
\end{Prop}
\begin{proof} Apply the facts that $h^{1,1}$ equals $10n$ \cite[Lemma IV.1.1]{MiES} and  that $h^{1,1}(X')-\rho_{\tr}(\pi')$ equals \[2n-2-\#\{\mbox{multiplicative fibers}\} - 2 \#\{\mbox{additive fibers}\} \]
(from e.g. \cite[Proposition 2.9]{Ext})    
to Lemma~\ref{dimLoc}.\end{proof}

\section{The lower bound} \label{under}
In this section we prove a lower bound for the dimension of $\NL_r$.

\begin{Thm}\label{locusdim}Let $r$ be an integer such that $2\leq r \leq 10n$. Let $L_r$ be the (constructible)  locus of Jacobian elliptic surfaces in $\M_n$ such that $\rho_{\tr}\geq r$ and the $j$-invariant is non-constant. Then
\[ \dim L_r = 10n-r. \]
\end{Thm}
\begin{proof}
Proposition~\ref{dimlem} implies that it suffices to  prove 
that there exists an elliptic surface without $II, III$ and $IV$ fibers, such that $\rho_{\tr}(\pi)=r$. {}From  Proposition~\ref{dimlem} it follows that such a surface lies on a component of $L_r$ of dimension $10n-r$. 

We start by choosing an elliptic surface $\pi_1:X_1\ra \Ps^1$ with four singular fibers, all multiplicative. By \cite[Proposition 2.9]{Ext} this fibration satisfies  $\rho_{\tr}(\pi_1)=10$. The existence of such surfaces has been established by Beauville \cite{Beaufour}; in particular there exist six such surfaces, up to isomorphism.

Let $\pi_{n}$ be a  cyclic base-change of degree $n$ of $\pi_{1}$ ramified at two points where the fibers are singular.

Since $\pi_{1}:X_{1} \ra \Ps^1$ satisfies $\rho_{\tr}(\pi_{1})=h^{1,1}(X_{1})$ (see e.g. \cite[Proposition 2.12]{Ext}), we obtain by \cite[Example 6.5]{Ext}
\[ \rho_{\tr}(\pi_{n})=h^{1,1}(X_{n})=10n, \]
which yields the claim in the case $r=10n$.

If $r<10n$, by the ``deformation of the $j$-map'' of $\pi_1$ (see \cite[Remark after Corollary~{3.5}]{MiRa} or combine Lemma~\ref{rambeh} with Lemma~\ref{ramred}) we can construct an elliptic surface $\pi : X \ra \Ps^1$ with
$2n+2+(10n-r)$ singular fibers, all multiplicative. By \cite[Proposition 2.9]{Ext} such a surface has
$\rho_{\tr}= r$. This finishes the proof.
\end{proof}

\begin{Cor}\label{locusdimCor} Let $r$ be an integer such that $2\leq r \leq 10n$. Then
\[ \dim \NL_r\geq 10n-r. \]
\end{Cor}

Another consequence of Theorem~\ref{locusdim} is the following: 

\begin{Cor}\label{K3cor} Let $MK3$ be the moduli space of  $K3$ surfaces. Let $2\leq r \leq 20$ be an integer.  Let $S_r$ be the locus in $MK3$ corresponding to $K3$ surfaces with $\rho(X)\geq r$. Then
\[ \dim S_r \geq 20-r. \]
\end{Cor}

\begin{proof} It is well-known that a Jacobian elliptic surface $\pi: X \ra \Ps^1$ with $p_g(X)=1$ is a $K3$ surface. 
Hence there is a morphism $\M_2 \rightarrow MK3$, which forgets the elliptic fibration. This morphism is  finite onto its image (see \cite{Sterk}). Let $C$ be a component of $L_r$ in $\M_2$ of dimension $20-r$. 
 The image of $C$ is contained in $S_r$ and is of dimension $20-r$.
\end{proof}

\begin{Rem} The surjectivity of the period map for (algebraic)
  $K3$ surfaces provides an alternative proof for the above
  result. Using the global Torelli theorem for $K3$ surfaces one
  obtains even equality. 
\end{Rem}

\section{Constant $j$-invariant}\label{constj}
We continue the study of $NL_r$ by considering the components of $\NL_r$ corresponding to
elliptic surfaces with constant $j$-invariant. In this section we
assume that $\pi: X \ra \Ps^1$ is an elliptic surface with precisely $2n$ fibers of type $I_0^*$.
All elliptic surfaces $\pi: X \ra \Ps^1$ with constant
$j$-invariant different from 0 or 1728, and $p_g(X)>0$ can be constructed in this way. The cases $j(\pi)=0$ and $j(\pi)=1728$ are discussed in the next section.

A Jacobian elliptic surface with $2n\; I_0^*$ fibers  is completely determined by the $2n$ points with an $I_0^*$ fiber and the $j$-invariant. Conversely, given a set $S$ of $2n$ points on $\Ps^1$ and a number $j_0\in \C-\{0,1728\}$ one can find a unique  elliptic surface (up to isomorphism) with $\pi: X \ra \Ps^1$ with $j(\pi)=j_0$ and $\Sing(\pi)=S$. (See Remark~\ref{RemHECon}.)
Hence the dimension of the (constructible) locus of all elliptic surface with $2n\;I_0^*$-fibers in $\M_n$ is $2n-2$, if $n\geq 2$.

\begin{Rem}\label{RemHECon} Let $\pi: X \ra \Ps^1$ be an elliptic surface with $2n$ fibers of type $I_0^*$. Then we associate with $\pi$ a hyperelliptic curve $\varphi: C \ra \Ps^1$ such that the ramification points of $\varphi$ are the points over which $\pi$ has a singular fiber. Let $E$ be an elliptic curve with the same $j$-invariant as the fibers of $\pi$. Then the minimal desingularization of $(C\times E)/\langle \iota \times [-1]\rangle$ is isomorphic to $X$. Conversely, given a hyperelliptic curve $C$ of genus $g$ the fibration induced by $(C\times E)/\langle \iota\times [-1]\rangle \ra C/\langle \iota \rangle\cong \Ps^1$ has constant $j$-invariant and $2g+2$ singular fibers of type $I_0^*$.   \end{Rem}

\begin{Rem}\label{RemSec} Consider  the elliptic surface $\pi: X \ra \Ps^1$ with 
\[ X=\widetilde{\frac{C\times E}{\langle \iota\times [-1]\rangle}}\] and  $\pi$ is induced by the projection  $C\times E \ra C$. 

Every section $s:\Ps^1\ra X$ comes from a morphism $\mu:C\ra E$ and $s$ maps a point $c\bmod \langle \iota \rangle$ to $(c,\mu(c)) \bmod \langle \iota\times [-1]\rangle$. Conversely a morphism $\mu$ defines a section if and only if $\mu$ maps the fixed points of $\iota$ to fixed points of $[-1]$. A constant morphism $\mu: C \ra \{P\}\subset E$ yields a section if and only if $P$ has order at most 2. This gives a contribution $(\Z/2\Z)^2$ to $\MW(\pi)$. Using \cite[Corollary VII.3.3]{MiES} one can show that $\MW(\pi)_{\tor}=(\Z/2\Z)^2$. If $\MW(\pi)\neq  (\Z/2\Z)^2$ then  a non-constant morphism $C\ra E$ exists with the above mentioned property.
\end{Rem}

\begin{Lem}\label{dimcal} Let $E$ be a curve of genus 1. Then the locus $L(E)$ corresponding to hyperelliptic curves $C$ admitting a non-constant morphism $C\ra E$ in $H_g$, the moduli space of hyperelliptic curves of genus $g$, has dimension $g-1$.
\end{Lem}

\begin{proof} {}From \cite[Lemma 1.1]{Scho} it follows that for any non-constant morphism $\psi: C \rightarrow E$, there exists an elliptic involution on $E$ induced by the hyperelliptic involution of $C$, i.e.,  such that the following diagram is commutative
\[ \begin{array}{ccc}
E &\leftarrow& C\\
\downarrow  &&\downarrow\\
\Ps^1 &\leftarrow &\Ps^1,\end{array}\]
where the vertical arrows are obtained by dividing out the (hyper)elliptic involution.

Fix $\lambda$ a Legendre parameter for $E$. Any non-constant morphism $f:\Ps^1 \rightarrow
\Ps^1$ gives rise to a  hyperelliptic curve $C=\widetilde{E \times_{\Ps^1} \Ps^1}$. The genus of $C$ is determined by $f$, i.e., $2g(C)+2$ equals the number of points with odd ramification index above $0,1,\lambda$ and $\infty$.

{}From this we obtain that $\dim L(E)$ equals the maximum over all $d$ of the dimension of the Hurwitz space corresponding to non-constant morphisms of degree $d$, such that above $0,1, \lambda, \infty$ there are precisely the $2g+2$ points with odd ramification index. By Corollary~\ref{ModN} this space has dimension $2\cdot 4-g-1+2g+2-2 \cdot 4-2=g-1$.
\end{proof}

\begin{Thm}\label{jconstThm} Let $n>1$. The locus $L$ in $\M_n$ of elliptic surfaces with $2n\; I_0^*$-fibers has dimension $h^{1,1}-\rho_{\tr}=2n-2=2p_g$. The locus $L_1$ of elliptic surfaces  with $2n \;I_0^*$ fibers and positive Mordell-Weil rank has dimension $p_g$. The locus $L_2$ of elliptic surfaces with $2n\; I_0^*$ fibers and Mordell-Weil rank at least 2 has dimension $p_g$ or $p_g-1$.
\end{Thm}

\begin{proof} A fiber of type $I_0^*$ has 4 components not intersecting the zero-section, so from the Shioda-Tate formula~\ref{ST} it follows that  $\rho_{\tr}=8n+2$. The first assertion follows  from the correspondence between $L$ and  sets of $2n$ distinct points in $\Ps^1$ together with a $j$-invariant $j_0\in \C$ as mentioned above.

For the second assertion, we note that by general theorems on the period map, the locus of elliptic surfaces with constant $j$-invariant and positive rank has dimension at most $p_g$. 
(One needs to exploit the well-known fact: for a cohomology class  $\xi \in H^2(X,\Z)$ lies in $H^{1,1}(X)$ if and only   $\xi\cdot \omega=0$, for every $\omega \in H^0(X,\Omega^2_X)$. Since $h^0(X,\Omega^2_X)=p_g$ this gives $p_g$ conditions on the image of the period map.) Hence $L_1$ has codimension at most $p_g$ in $L$.

If $\MW(\pi)$ is strictly bigger than $(\Z/2\Z)^2$ then by Remark~\ref{RemSec} there is a
non-constant morphism $C \rightarrow E$, with $C$ and $E$ as in Remark~\ref{RemHECon}. Hence for a fixed $j_0\in \C$, the locus of elliptic surfaces with $j(\pi)=j_0$ and
$\MW(\pi)$ infinite has by Lemma~\ref{dimcal} dimension at most $g(C)-1$. Hence $L_1$ has dimension at most $g(C)=p_g(X)$.

Suppose the the fixed $j_0$ corresponds to a curve with complex multiplication. Since the
Mordell-Weil group of $\pi$ modulo torsion is a free $\End(E)$-module, it follows that $\rank MW(\pi)$ is even, so $L_2$ has dimension at least $p_g-1$.
\end{proof}

\section{$j$-invariant 0 or  1728} \label{jnul}
In this section we will prove that  $\dim \NL_r-(10n-r)$ can be arbitrarily large.

\begin{Prop}\label{ExisPropA} Let $n\geq 2$. Fix an integer $k$ such that $6n/5\leq k\leq 6n$. Then there exists an elliptic surface $\pi: X \ra \Ps^1$ with $j(\pi)=0$, $p_g(X)=n-1$ and $k$ singular fibers. Moreover, the locus of elliptic surfaces with $j(\pi')=0$ and $C(\pi')=C(\pi)$ has dimension $k-3$ in $\M_n$. If $m$ is an integer such that $m>6n$ or $m<6n/5$ then there exists no elliptic surface with $j(\pi')=0$ and $m$ singular fibers.\end{Prop}

\begin{Prop}\label{ExisPropB}Let $n\geq 2$.
Fix an integer $k$  such that  $4n/3\leq k\leq 4n$. Then there exists an elliptic surface $\pi: X \ra \Ps^1$ with $j(\pi)=1728$, $p_g(X)=n-1$ and $k$ singular fibers. Moreover, the locus of elliptic surfaces with $j(\pi')=1728$ and $C(\pi')=C(\pi)$ has dimension $k-3$ in $\M_n$. If $m$ is an integer such that $m>4n$ or $m<4n/3$ then there exists no elliptic surface with $j(\pi')=1728$ and $m$ singular fibers.\end{Prop}

\begin{proof}[Proof of Propositions~\ref{ExisPropA} and~\ref{ExisPropB}] Without loss of generality we may assume that all elliptic surfaces under consideration have a smooth fiber over $\infty$.

An elliptic surface with $k$ singular fibers, $p_g(X)=n-1$ and $j(\pi)=0$ exists if and only if there exists a polynomial $f$ of degree $6n$ with $k$ distinct zeroes, and every zero has multiplicity at most 5. We can associate with any such  polynomial $f(t)$ an elliptic surface with Weierstrass equation $y^2=x^3+f(t)$, and vice-versa, an elliptic surface with $j$-invariant $0$ gives rise to a Weierstrass equation of the above form.

Hence an elliptic surface with $k$ singular fibers exists if and only if
$6n/5\leq k\leq 6n$.
Modulo the action of $\Aut(\Ps^1)$ we obtain a $k-3$ dimensional locus in $\M_n$.

The case of $j(\pi)=1728$ is similar except for the fact that the polynomial $g(t)$ is of degree $4n$, and the highest possible multiplicity is 3. The associated surfaces is then given by $y^2=x^3+g(t)x$. \end{proof}

\begin{Prop}\label{jnulprop} Let $n\geq 2$. Let $r\leq 1 +\frac{24}{5} n$ be a positive integer. Then the locus of elliptic surfaces with $j$-invariant $0$ and $\rho_{\tr}(X)$ at least $2r$ has dimension
\[ 6n-r-2 \]
\end{Prop}

\begin{proof} If $j(\pi)$ is constant and $\pi$ has $k$ singular fibers then the number of components of singular fibers not intersecting the zero-section equals $12n-2k$ (see e.g. \cite[Proposition 2.12]{Ext}). Hence $\rho_{\tr}(\pi)=2+12n-2k$. {}From this it follows that $\rho_{\tr}(\pi)\geq 2r$ if and only if $k\leq 6n-r+1$. We want to apply Proposition~\ref{ExisPropA} for $k=6n-r+1$. The condition on $k$ is equivalent  to the above assumption on $r$. Then Proposition~\ref{ExisPropA} implies that the dimension of the locus is $k-3=6n-r-2$.
\end{proof}

\begin{Rem} A similar result holds in the case that $j(\pi)=1728$. In that case one should take $r\leq \frac{14}{3} n+1$.\end{Rem}

\begin{Rem} All loci $L$ described in the Sections~\ref{under} and~\ref{constj} satisfied $\dim L + \rho(X) \leq 10 n$, for an $X$ corresponding to a generic point of $L$.
In Proposition~\ref{jnulprop} one can choose $r=1+ 4n +\lfloor 4n/5 \rfloor$,
with $\lfloor \alpha \rfloor$ denoting the largest integer, not larger then $\alpha$. One obtains
\[ \dim L + \rho(X) =6n-r-2+2r=10n +\left\lfloor \frac{4}{5} n\right\rfloor - 1 \]
The excess term $\lfloor 4n/5 \rfloor-1$ can be arbitrarily large.
\end{Rem}

\begin{Cor} Suppose $n\in \{ 2,3,4,5 \}$. Then $\dim \NL_{10n} = n-2$.
\end{Cor}

\begin{proof} {}From the Infinitesimal Torelli theorem for Jacobian elliptic surfaces \cite[Corollary 4.3]{Ext} 
it follows that \[ \{ [\pi:X \ra \Ps^1] \in \NL_{10n} \mid \rho_{\tr}(\pi)<10n \mbox{ or } j(\pi)\mbox{ not constant}\} \]
  is a discrete set. If $j(\pi)\in \C - \{ 0,1728 \}$ then
  $\rho_{\tr}(\pi)=8n+2<10n$, hence we only have to consider elliptic
  surfaces with $\rho_{\tr}(\pi)=10n$ and constant $j$-invariant 0 or
  1728. Since $1+\lfloor 24n/5\rfloor=5n$ for the $n$ under consideration, we may apply Proposition~\ref{jnulprop} with $r=5n$. This yields  $\dim \NL_{10n}=n-2$.
\end{proof}

\begin{Rem} {}From this Corollary we deduce that for  $n \in \{3,4,5\}$, there exist positive dimensional loci $L \subset \M_n$, such that any surface $X$ corresponding to a point in $L$ satisfies $\rho(X)=h^{1,1}(X)$. The image of the period map restricted to $L$ has  discrete image. This contradicts several  Torelli-type properties (see also \cite[Theorem 4.8]{Ext}).
\end{Rem}

\section{Upper bound}\label{upper}

As in \cite{Cox}, we study the Noether-Lefschetz loci using the
identification of $H^{1,1}, H^{2,0}$ and $H^{0,2}$ with several
graded pieces of a Jacobian ring $R$. We choose to give a more
algebraic presentation than in \cite{Cox}.

To be precise, given a Weierstrass minimal equation $F=0$ for $\pi: X \ra
\Ps^1$ we can construct a (singular) hypersurface $Y$ in the weighted
projective space $\Ps:=\Ps{} (1,1,2n,3n)$ (with projective
coordinates $x$, $y$, $z$, $w$ of weight 1, 1, $2n$, $3n$ resp.) given by:
\[ 0=-w^2+z^3+P(x,y)z+Q(x,y)=:F \]
with  $n=p_g(X)+1$, $\deg(P)=4n$ and  $\deg(Q)=6n$. 
Here $X$ and $Y$ are birational; $Y$ is obtained from $X$ by contracting the
zero-section and all fiber components not intersecting the zero-section. 

Let $A:=\C[x,y,z,w]$ with weights $1,1,2n,3n$. Let $B=\C[x,y]\subset A$. The construction $(\pi:X\ra \Ps^1, \sigma_0:\Ps^1\ra X) \mapsto Y$ gives a nice description of the moduli space $\M_n$. (See the proof of Theorem~\ref{upperthm}.)

Let $J\subset A$ be the ideal generated by the partial derivatives of $F$. 
The Jacobi ring $R$ is the quotient ring $A/J$. 
 It is well known (see \cite{Cox}, \cite{CD}, \cite{Dol}, \cite{Ste}) that if all the fibers of $\pi$ are irreducible then $Y$ is quasi-smooth, i.e., the cone $(F=0)\subset \A^4$ is smooth outside the origin. 

Assume for the moment that $\pi$ satisfies this assumption, i.e., $Y$ is quasi-smooth. 
Then the classical Griffiths-Steenbrink theorem, applied to our case, states that we have isomorphisms
\[ H^{2,0}(Y) \cong R_{n-2}, \; H^{1,1}(Y)_{\prim} \cong R_{7n-2}, \; H^{0,2}(Y) \cong R_{13n-2}.\]
Here, we adopt the convention that for a graded ring $R'$ we denote by $R'_d$ all elements of degree $d$ and for a variety $Y'\subset \Ps$ we denote by $H^{1,1}(Y')_{\prim}=\Ima (H^2(\Ps,\C) \ra H^{1,1}(Y'))^{\perp}$.

In the case that $\pi$ has reducible fibers the situation is very
similar.
This follows from a special case of a recent result of Steenbrink \cite{Stermk}.  Note that in this case $Y$ is not quasi-smooth. 

\begin{Thm}[{Steenbrink \cite{Stermk}}]\label{SteThm} Let $Y'\subset \Ps$ be a surface of degree $6n$, whose only singularities outside $\Ps_{\sing}$ are rational double points and which is transverse to   $\Ps_{\sing}$.
Then there is  a natural isomorphism
$H^{2,0}(Y') \cong R'_{n-2}$ and an injective map
\[ H^{1,1}(Y')_{\prim} \ra R'_{7n-2}.\]  
\end{Thm}

\begin{Lem}\label{redlem} We have
\[ H^{1,1}(Y)_{\prim} \cong H^{1,1}(X)/(T(\pi)\otimes \C).\]
In particular, $\dim H^{1,1}(Y)_{\prim} = 10n-\rho_{\tr}$.
\end{Lem}
\begin{proof}
The 
  isomorphism follows from the fact that $\varphi: X\ra Y$ is a resolution of singularities, $\varphi$ contracts the zero-section and all fiber components not intersecting the zero-section
and the fact that a general hyperplane section $H\cap Y$ is a fiber of $\pi$.
\end{proof}

\begin{Cor}\label{GrifCor} There is a natural isomorphism
$H^{2,0}(X) \cong R_{n-2}$ and an injective map
\[ H^{1,1}(X)/(T(\pi)\otimes \C ) \ra R_{7n-2}.\]
\end{Cor}

\begin{proof} The existence of the two linear maps follows from Theorem~\ref{SteThm} and Lemma~\ref{redlem}.
\end{proof}

Next, we prove some elementary technical results. For a
polynomial $P$, we use a subscript (like $P_x$) to indicate the
derivative with respect to the variable in the subscript.

\begin{Lem}\label{derlem} Let $\pi: X \ra \Ps^1$ be the elliptic surface associated with
$w^2=z^3+Pz+Q$, with $P\in \C[x,y]_{4n}, Q\in \C[x,y]_{6n}$. Then
$P_xQ_y-P_yQ_x=0$ if and only if $j(\pi)$ is constant.\end{Lem}

\begin{proof}[Proof (see \cite{Cox})] Using the Euler relation for weighted homogeneous polynomials one easily obtains that the partial derivative to $x$ or to $y$ of $j(\pi)=1728\cdot 4
P^3/(4P^3+27 Q^2)$ is identically zero if and only if
$(P_xQ_y-P_yQ_x) P Q=0$. If $PQ$ is zero then also $P_xQ_y-P_yQ_x$ equals zero, which
gives the lemma.
\end{proof}

\begin{Lem} \label{Jacobi-rank}Fix a positive integer $n$. Let $F\in
  A=\C[x,y,z,w]$ be a weighted homogeneous polynomial of degree
  $6n$. Suppose  that the variety in $\Ps(1,1,2n,3n)$ defined by $F=0$
  is birational to an elliptic surface $\pi: X \ra \Ps^1$, with $\pi$
  induced by $[x,y,z,w]\mapsto [x,y]$ and $F=0$ is  a Weierstrass minimal
  equation. Let $J$ be the Jacobi-ideal of $F$. Let $\widetilde{J}\subset A$ be the $B$-submodule generated by $J_{\leq 6n}$.  If $j(\pi)$ is not constant then $\widetilde{J}$ is a free $B$-module of rank 7, otherwise it is a $B$-module of rank 6.
\end{Lem}

\begin{proof}
After applying an automorphism of $\Ps$ we may assume that $F=-w^2+z^3+Pz+Q$,
with $P\in B_{4n}$ and $Q\in B_{6n}$.  Then we have the following set of generators for $\widetilde{J}$ as a $B$-module
\[ w^2,wz, w, 3z^3+Pz, 3z^2+P, P_x z+Q_x, P_yz+Q_y.\]
By degree considerations, we obtain that the first six generators generate a free $B$-module of rank 6. It suffices to prove that $\widetilde{J}$ has rank 6 if and only if $j(\pi)$ is constant.

Consider the elements
$\alpha:=Q_y(P_xz+Q_x)-Q_x(P_yz+Q_y)=(P_yQ_x-P_xQ_y)z $ and
$\beta:=P_y(P_xz+Q_x) - P_x(P_yz+Q_y)=   P_yQ_x-P_xQ_y $.

Suppose $j(\pi)$ is constant then we obtain by Lemma~\ref{derlem}
the relation $P_y(P_xz+Q_x) -P_x(P_yz+Q_y)=0$, proving that the rank of $\widetilde{J}$ is 6.

Suppose $j(\pi)$ is not constant. Then Lemma~\ref{derlem} implies
that $\alpha$ and $\beta$ are non-zero, hence are independent. By degree considerations we obtain that
\[ w^2, wz, w,3z^3+Pz,3z^2+P,\alpha, \beta\]
generate a free submodule of $\widetilde{J}$ of rank 7. This finishes the proof.
\end{proof}

\begin{Rem} Let $\widetilde{A}$ be the $B$-module generated by $A_{\leq 6n}$. Then $\widetilde{A}$ has  rank 7. So $\widetilde{J}$ has the same rank as $\widetilde{A}$ if and only if $j(\pi)$ is not constant.\end{Rem}

\begin{Prop}\label{Jacobi-codim} Let $\pi: X \ra \Ps^1$ be an elliptic surface such that $j(\pi)$ is not constant. Let $Y\subset \Ps(1,1,2n,3n)$ be the associated surface in weighted projective space. Let $V\subset A_{6n}$ be a vector space containing the degree $6n$-part of the Jacobi-ideal of $Y$. Let $V_k$ be the image of the multiplication map $V\otimes B_k \ra \widetilde{A}_{6n+k}$. Then for all $k\geq 0$
\[ \codim_{\widetilde{A}_{6n+k}} V_k \leq \codim_{\widetilde{A}_{6n}} V.\] 
\end{Prop}

\begin{proof} 
Since $\widetilde{J}$ and $\widetilde{A}$ are free of rank 7 and generated by elements of degree at most $6n$, we obtain \begin{eqnarray*}\dim \widetilde{J}_{6n+k}&=&\dim J_{6n}+7k\\ \dim \widetilde{A}_{6n+k}&=&\dim A_{6n}+7k.\end{eqnarray*}

Let $\widetilde{V} :=\oplus_{k\geq 0} V_k$.
Observe that
\[ \dim \widetilde{J}_{6n+k}\leq \dim V_k\leq \dim \widetilde{A}_{6n+k}.\]
Using that the Hilbert function of $\widetilde{V}$ is a polynomial, these inequalities imply that $\dim V_k=c+7k$, for $k\gg0$.
Let $d_k:=\dim \widetilde{V}_{k+1}-\dim \widetilde{V}_k$. Since $V$ is a torsion-free $B$-module generated in degree 0 we have that $d_k$ is a decreasing function for $k\geq0$. {}From this we obtain $d_k\geq 7$ for $k\geq 0$.
Hence $\dim V_{k}\geq \dim V+7k$. Recall that $\dim \widetilde{A}_{6n+k}=\dim
\widetilde{A}_{6n}+7k$, we obtain that $\codim_{\widetilde{A}_{6n}} V\geq \codim_{\widetilde{A}_{6n+k}} V_{k}$, which finishes the proof.
\end{proof}

Theorem~\ref{upperthm} together with the results in the previous sections will provide a proof for Theorem~\ref{MainThmNL}.

\begin{Thm}\label{upperthm}  Let $2\leq r \leq 10n$. Let $U\subset \M_n$ be the locus of elliptic surfaces with non-constant $j$-invariant. Then $\dim \NL_r\cap \; U$ is at most $10n-r$.
\end{Thm}

\begin{proof}
We prove the theorem by descending induction. Assume that it is true for all $r'$, $r<r'\leq 10n+1$. 

Define $C$ by  $\NL_r=\NL_{r+1} \coprod C$. Then by induction we have that the dimension of $NL_{r+1}$ is at most $10n-r-1$. Hence it suffices to prove that $\dim C\leq 10n-r$.
Let $\pi:X\ra \Ps^1$ correspond to a point $p$ in $C\cap U$. We want to calculate the
dimension of the tangent space of $\NL_r$ at $p$.
Let $Y\subset \Ps$ be the
corresponding surface in the weighted projective space $\Ps$. Write $s=\rank \MW(\pi)$.

 The moduli space $\M_n$ can be obtained in the following way: Let
\[U:=\left\{f\in A_{6n} \colon \begin{array}{c} f=0 \mbox{ is birational to an elliptic
  surface } \\ \mbox{and } f \mbox{ is Weierstrass minimal. }\end{array}\right\} \]
then $U/\Aut(\Ps)=\mathcal{M}_n$ (see \cite{MiMS}).
Let $L\subset A_{6n}$ be the pre-image of a component containing $\pi: X \ra \Ps^1$ of
$\NL_r\subset \M_n$ .

Using Miranda's construction of $\M_n$ it follows that the codimension of $L/\Aut(\Ps)$ in $\M_n$ equals $\codim_{A_{6n}}L$. {}From this it follows that it suffices to show that $L$ has codimension
at least $r-2$ in $A_{6n}$. Let $\T\subset A_{6n}$ be the tangent space to $L$ at $Y$, considered as a point in  $A_{6n}$.

Consider the multiplication map
\[ \varphi : \T \otimes A_{n-2} \ra A_{7n-2}.\]
Let $\psi$ be the composition of $\varphi$ with the projection onto $R_{7n-2}$.
Using Corollary~\ref{GrifCor} 
we obtain that $\psi$ corresponds to the map $\mathcal{T} \otimes H^{2,0}\ra H^{1,1}_{\prim}$ induced by the period map. Hence the image of $\varphi$ is contained in the subspace $W\subset A_{7n-2}$ that is the pre-image of $H^{1,1}_{\prim}\hookrightarrow R_{7n-2}$ (using Corollary~\ref{GrifCor}).
We have that $J_{7n-2}$ is contained in $W$. {}From Lemma~\ref{redlem} it follows that 
$\dim R_{7n-2}=10n-2$. Since $\codim_{R_{7n-2}} H^{1,1}_{\prim}=\rho_{\tr}-2$ it follows that the difference  $\dim W - \dim J_{7n-2}$ equals $10n-\rho_{\tr}$.  {}From this we obtain $\codim_{A_{7n-2}} W=\rho_{\tr}-2$.

Consider the sheaf $\mathcal{H}^{1,1}  : Y'\mapsto H^{1,1}(Y',\C)_{\prim}$ on $U\subset A_{6n}$, which is a subsheaf of the constant sheaf $Y'\mapsto H^2(Y',\Z)_{\prim}\otimes\C$ on $U$.

Let  $O$ be the Zariski-constructible set of $Y'\in L$ such that $\rho_{\tr}(\pi')=r-s$ and the rank of $ MW(\pi')$ equals $s$, where we use the fibration $\pi'$ induced by projection from the singular point of $Y'$. Then on $O$  there is a constant subsheaf $\mathcal{NS}$ of $\mathcal{H}^{1,1}_{\prim}$, given by $Y' \mapsto  NS(Y')_{\prim} \otimes \C$. The stalks of this sheaf are isomorphic to $\C^s$. 
{}From this it follows that the image of the contraction
$\T\otimes H^{2,0} \ra H^{1,1}$ is orthogonal (for the cup-product) to the stalk of $\mathcal{NS}$ at $p$. In particular, it has codimension at least $s$.
Using Theorem~\ref{SteThm} it follows  that   $\codim_W \T_{n-2}\geq s$, hence \[\codim_{A_{7n-2}} \T_{n-2} \geq r-2,\] where $\T_{n-2}$ denotes the image of $\varphi$ restricted to $\T\otimes A_{n-2}$.

Hence it suffices to show that 
\[\codim_{A_{6n}} \T \geq \codim_{A_{7n-2}} \T_{n-2}.\]
Since $L$ is stable under $\Aut(\Ps)$ , its tangent space $\T$ contains the subspace induced by the Lie algebra of $\Aut(\Ps)$. One can show that this subspace is $J_{6n}$. Hence we can apply Lemma~\ref{Jacobi-codim} with $V=\T$ and $k=n-2$. Using $\widetilde{A}_{7n-2}={A}_{7n-2}$ we obtain the desired inequality.   
\end{proof}
\section{Concluding remarks}\label{concl}

\begin{Rem} The argument used in the proof of Theorem~\ref{upperthm}
  cannot work for elliptic surfaces with constant $j$-invariant. First
  of all, in this case  Lemma~\ref{Jacobi-rank} gives  $d_k\geq
  6$, which implies only $\codim_{A_{6n}} V\geq
  \codim_{A_{7n-2}} V_{n-2} -(n-2)$. Moreover, it is not  hard to give
  a linear subspace $V\subset A_{6n}$ such that $J_{6n}\subset V$,
  $V\neq J_{6n}$ and $\codim_{A_{7n-2}} V_{n-2} >\codim_{A_{6n}} V$. One
  needs to show that such spaces do not occur as the tangent space to
  a component of $\NL_r$, different from the 
  components already described in Section~\ref{jnul}. By the results of Section~\ref{constj} we know that such a $V$ would have a large codimension in $A_{6n}$, but these results are not sufficient to prove the theorem in the case of constant $j$-invariant.
\end{Rem}

\begin{Rem}
There is still an interesting open issue. In the theory of Noether-Lefschetz loci there is the notion of special components and of general components. Special components are the components of $\NL_3$ with codimension in $\NL_2$ less then $p_g$. In the case of elliptic surfaces there is only one special component (see \cite{Cox}). For higher Noether-Lefschetz loci, one can define the special components as the components in $\NL_r$ with codimension less then $(r-2)p_g$. Then one finds infinitely many special components. One can also define special components as the components of $\NL_r$ such that the maximal codimension in $\NL_{r-1}$ is less then $p_g$. By base-changing families of elliptic $K3$ surfaces we can  find again infinitely many special components, even when we fix the component of $\NL_{r-1}$ in which these components are contained.
\end{Rem}
\begin{Rem} Suppose $\M$ is a  moduli space for some class of smooth
  surfaces. We
would like to
obtain $\codim_\M \NL_r \geq r-\rho_{\gen}$, where $\rho_{\gen}$ stands
  for the generic Picard number and $\NL_r=\{X \in \M\mid \rho(X)\geq r\}$.

To give a proof similar to the proof of 
  Theorem~\ref{upperthm} it suffices to assume the following conditions:

\begin{itemize}
\item Griffiths-Steenbrink holds for the moduli problem. I.e.,
there exists a threefold $X$, such that for all points $p\in
  \M$ there exists  a surface $Y_p\subset X$,
  satisfying the conditions of
   \cite{Stermk}. Moreover, if $\widetilde{Y}_p$ is the
  normalization of $Y$, then $[\widetilde{Y}_p]\in \M$ is the point $p$.
  
\item All surfaces are linearly equivalent (as divisors on $X$), i.e., fix a point $p\in \M$. Let $X$ and $Y_p$ as above, then there is a dense open $U\subset H^0(X,\str_X(Y_p))$ and a
  surjective morphism   $U \ra \M$, sending a divisor $Y'$ to the
  class of its minimal desingularization.
\item The following  multiplication conditions hold. Let $K$ be the kernel of
  $\psi_2$. Let $K(m)$ be the image of $K^{\otimes m}$ in
  $H^0(X,K_X^{\otimes m}(2mY))$. Then for all $m\geq 2$ we have
\begin{eqnarray*} && \dim K(m)-\dim K(m-1)\geq\\ & &\geq \dim H^0(X,K_X^{\otimes m}(2mY))-\dim
  H^0(X,K_X^{\otimes (m-1)}(2(m-1)Y)).\end{eqnarray*} 
\end{itemize}
\end{Rem}

\begin{Rem} In \cite{GreenF}, the following statement is proven. Let
  $d>3$ be an integer, let
  $U\subset \C[x,y,z,w]_d$ be the set of homogeneous polynomials $F$
  such that $F=0$ defines a smooth surface. Let $\NL\subset U$ be the
  locus of surfaces with Picard number at least 2. Then $\codim_U \NL=d-3$.

The strategy used in the proof is very similar to the strategy used in
the proof of Theorem~\ref{upperthm}. However, in this case the
strategy does not seem to work for larger Picard number. If one
applies a reasoning as in the proof of Theorem~\ref{upperthm} one obtains
that $\codim_U \NL_r \geq d-1-r$. Griffiths and Harris \cite[page
208]{GrHaCp} {\em conjecture} that for $3\leq r \leq d$ we have
\[ \codim_U \NL_r=(r-1)(d-3)-\left( \begin{array}{c} r-3 \\2\end{array}
  \right).\] 
and they claim that it is easy to prove that we can replace the
equality sign by a less or equal sign.

There is still a gap between these two bounds for $\codim_U \NL_r$.
\end{Rem}

\end{document}